\documentclass{amsart}

\usepackage{amscd}

\theoremstyle{plain}
\newtheorem{theorem}{Theorem}[section]

\newtheorem{lemma}[theorem]{Lemma}

\theoremstyle{definition}
\newtheorem{definition}[theorem]{Definition}

\theoremstyle{remark}
\newtheorem{remark}[theorem]{Remark}

\newcommand{\N}{\mathbb{N}}                   
\newcommand{\R}{\mathbb{R}}                   
\newcommand{\C}{\mathbb{C}}                   
\newcommand{\K}{\mathbb{K}}                   
\newcommand{\Oh}{\widehat{\mathcal O}}        
\newcommand{\una}{\underline{a}}              
\newcommand{\vp}{\varphi}
\newcommand{\vph}{\hat{\varphi}}
\newcommand{\mi}{\mathfrak m}                 
\newcommand{\mih}{\widehat{\mi}}              
\newcommand{\Ni}{\mathfrak N}                 
\newcommand{\Runa}{{\mathcal R}_{\una}}       
\newcommand{\Ra}{{\mathcal R}_a}
\newcommand{\ua}{\underline{a}}
\newcommand{\ux}{\underline{x}}
\newcommand{\uvp}{\underline{\varphi}}
\newcommand{\Ms}{{M^s_{\varphi}}}
\newcommand{\Msr}{{\mathring{M}^s_{\varphi}}}
\newcommand{\Es}{{E^s_{\varphi}}}
\newcommand{\Esr}{{\mathring{E}^s_{\varphi}}}

\newcommand{\Fb}{{\mathcal F}_b}              
\newcommand{\supp}{\mathrm{supp}\,}           

\newcommand{\ad}{\mathrm{ad}\,}

\newcommand{\Hom}{\mathrm{Hom}}
\newcommand{\Imm}{\mathrm{Im}\,}
\newcommand{\Ker}{\mathrm{Ker}\,}
\newcommand{\id}{\mathrm{id}}
\newcommand{\rk}{\mathrm{rk}\,}
\newcommand{\NR}{\mathrm{NR}}

\newcommand{\cO}{\mathcal O}
\newcommand{\cC}{\mathcal C}
\newcommand{\al}{\alpha}
\newcommand{\be}{\beta}
\newcommand{\lbr}{{[\![}}
\newcommand{\rbr}{{]\!]}}

\numberwithin{equation}{section}

\begin{document}

\title{Uniform linear bound in Chevalley's lemma}
\author{J. Adamus, E. Bierstone and P.D. Milman}
\thanks{Research partially supported by NSERC Postdoctoral Fellowship
PDF 267954-2003 (Adamus), and NSERC Discovery Grants OGP 0009070 (Bierstone),
OGP 0008949 (Milman)}

\address{J. Adamus, Institute of Mathematics of the Polish Academy of 
Sciences, 00-956 Warszawa 10, Sniadeckich 8, P.O. Box 21, Poland}
\email{adamus@impan.gov.pl}

\address{E. Bierstone, Department of Mathematics, University of Toronto, 
Toronto, Ontario, Canada M5S 3G3}
\email{bierston@math.toronto.edu}

\address{P.D. Milman, Department of Mathematics, University of Toronto,
Toronto, Ontario, Canada M5S 3G3}
\email{milman@math.toronto.edu}

\keywords{Chevalley function, regular mapping, Nash subanalytic set}

\begin{abstract}
We obtain a uniform linear bound for the Chevalley function at a point in
the source of an analytic mapping that is regular in the sense of
Gabrielov. There is a version of
Chevalley's lemma also along a fibre, or at a point of the image of a proper
analytic mapping. We get a uniform linear bound for the Chevalley
function for a closed Nash (or formally Nash) subanalytic set.
\end{abstract}
\maketitle
\setcounter{tocdepth}{1}
\tableofcontents

\section{Introduction}
\label{sec:intro}

Let $\vp : M\to N$ denote an analytic mapping of analytic manifolds
(over $\K = \R$ or $\C$). Let $a\in M$. Let $\vp^*_a : \cO_{\vp(a)}
\to \cO_a$ or $\vph^*_a : \Oh_{\vp(a)} \to \Oh_a$ denote the induced
homorphisms of analytic local rings or their completions, 
respectively. (We write $\cO_a = \cO_{M,a}$, and $\mi_a$ (or $\mih_a$) $=$ 
maximal ideal of $\cO_a$ (or $\Oh_a$).) According to
Chevalley's lemma (1943), there is an increasing function $l: \N \to \N$
(where $\N$ denotes the nonnegative integers) such that
$$
\vph^*_a(\Oh_{\vp(a)}) \cap \mih_a^{l(k)+1}\ \subset\ \vph^*_a(\mih_{\vp(a)}
^{k+1})\ ;
$$
i.e., if $F \in \Oh_{\vp(a)}$ and $\vph^*_a(F)$ vanishes to order
$l(k)$, then $F$ vanishes to order $k$, modulo an element 
of $\Ker \vph^*_a$ (\cite{Chev}; cf. Lemma \ref{lem:Chevalley} below). 
Let $l_{\vp^*}(a,k)$ denote the least $l(k)$ satisfying Chevalley's lemma.
We call $l_{\vp^*}(a,k)$ the \emph{Chevalley function} of $\vph^*_a$.

Let $x = (x_1, \ldots, x_m)$ and $y = (y_1, \ldots, y_n)$ denote local
coordinate systems for $M$ and $N$ at $a$ and $\vp(a)$, respectively.
The local rings $\cO_a$ or $\Oh_a$ can be identified with the rings
of convergent or formal power series $\K\{x\} = \K\{x_1, \ldots, x_m\}$
or $\K\lbr x\rbr = \K\lbr x_1, \ldots, x_m\rbr$, respectively. In the local
coordinates, write $\vp(x) = (\vp_1(x), \ldots, \vp_n(x))$. Then $\Ker 
\vph^*_a$ is the \emph{ideal of formal relations} $\{F(y) \in 
\K\lbr y\rbr:\ F(\vp_1(x), \ldots, \vp_n(x)) = 0\}$ (and $\Ker \vp^*_a$ is
the analogous \emph{ideal of analytic relations}). Chevalley's lemma is an
analogue for such nonlinear relations of the Artin-Rees lemma. (See Remark
\ref{artinrees}.)

Let $r^1_a(\varphi)$ denote the generic rank of $\varphi$ near $a$, and set
\[
r^2_a(\varphi)\ :=\ \dim\frac{\Oh_{\varphi(a)}}{\Ker\vph^*_a}\ ,\qquad 
r^3_a(\varphi)\ :=\ \dim\frac{\mathcal{O}_{\varphi(a)}}{\Ker{\varphi}^*_a}
\]
(where $\dim$ denotes the Krull dimension). Then $r^1_a(\vp)\leq r^2_a(\vp)
\leq r^3_a(\vp)$. Gabrielov proved that if $r^1_a(\vp) = r^2_a(\vp)$, then
$r^2_a(\vp) = r^3_a(\vp)$ \cite{Ga}; i.e., if there are enough formal relations,
then the ideal of formal relations is generated by convergent relations. 
The mapping $\vp$ is called \emph{regular at} $a$
if $r^1_a(\vp)=r^3_a(\vp)$. We say
that $\vp$ is \emph{regular} if it is regular at every point of $M$. Izumi
\cite{Iz2} proved that $\varphi$ is regular at $a$ if and only if the 
Chevalley function of $\vph^*_a$ has a \emph{linear (upper) bound};
i.e., there exist $\al, \be \in \N$
such that 
$$
l_{\vp^*}(a,k)\ \leq\ \al k + \be\ ,
$$
for all $k \in \N$. On the other hand, Bierstone and Milman \cite{BM3} proved
that, if $\vp$ is regular, then $l_{\vp^*}(a,k)$ has a \emph{uniform bound};
i.e., for every compact $L \subset M$, there exists $l_L:\ \N \to \N$
such that
$$
l_{\vp^*}(a,k)\ \leq\ l_L(k)\ ,
$$
for all $a \in L$ and $k \in \N$. In this article, we prove that the 
Chevalley function associated to a regular mapping has a 
\emph{uniform linear bound}:

\begin{theorem}
\label{thm:main1}
Suppose that $\vp$ is regular. Then, for every compact $L \subset M$, there
exist $\al_L, \be_L \in \N$ such that
$$
l_{\vp^*}(a,k)\ \leq\ \al_L  k + \be_L\ ,
$$
for all $a \in L$ and $k \in \N$.
\end{theorem}

Chevalley's lemma can be used also to compare two notions of order
of vanishing of a real-analytic function at a point of a subanalytic set.
Let $X$ denote a closed subanalytic subset of $\R^n$. Let $b \in X$ and
let $\Fb(X) \subset \R\lbr y-b\rbr$ denote the formal local ideal of 
$X$ at $b$. (See Lemma~\ref{lem:formal-closure}.) For all $F \in \Oh_b = 
\R\lbr y-b\rbr$, we define
\begin{equation}
\label{eqn:ord}
\begin{aligned}
\mu_{X,b}(F) &:= \max\{l\in\N:\  |T^l_bF(y)|\leq\mathrm{ const\,}|y-b|^l,
\ y\in X\}\ ,\\
\nu_{X,b}(F) &:= \max\{l\in\N:\ F\in\mih_b^l+\Fb(X)\}\ ,
\end{aligned}
\end{equation}
where $T^l_bF(y)$ denotes the Taylor polynomial of order $l$ of $F$ at $b$.
Then there exists $l:\ \N \to \N$ such that, for all $k \in \N$, 
if $F\in\Oh_b$ and $\mu_{X,b}(F) > l(k)$, then $\nu_{X,b}(F) > k$.
(See Section 3.) 
For each $k$, let
$l_X(b,k)$ denote the least such $l(k)$. We call $l_X(b,k)$ the 
\emph{Chevalley function} of $X$ at $b$.

\begin{theorem}
\label{thm:main2}
Suppose that $X$ is a Nash (or formally Nash) subanalytic subset of $\R^n$.
Then the Chevalley function of $X$ has a \emph{uniform linear bound}; 
i.e., for every
compact $K \subset X$, there exists $\al_K, \be_K \in \N$ such that 
$$
l_X(b,k)\ \leq\ \al_K k + \be_K\ ,
$$
for all $b \in K$ and $k \in \N$.
\end{theorem}

Theorems \ref{thm:main1} and \ref{thm:main2} are the main new results in this
article. They answer questions raised in \cite[1.28]{BMAnn}.

The closed \emph{Nash subanalytic} subsets $X$ of $\R^n$ are the images
of regular proper real-analytic mappings $\vp:\ M \to \R^n$. In particular,
a closed semianalytic set is Nash.
A closed subanalytic subset 
$X$ of $\R^n$ is \emph{formally Nash} if, for every $b \in X$, there is
a closed Nash subanalytic subset $Y$ of $X$ such that $\Fb(X) = \Fb(Y)$
\cite{BMAnn}. Unlike the situation of Theorem \ref{thm:main1}, the
converse of Theorem \ref{thm:main2} is false \cite[Example 12.8]{BMAnn}.

The main theorem of \cite{BMAnn} (Theorem 1.13) asserts that, if $X$ is
a closed subanalytic subset of $\R^n$, then the existence of a uniform
bound for $l_X(b,k)$ is
equivalent to several other natural analytic and algebro-geometric conditions; for
example, semicoherence \cite[Definition 1.2]{BMAnn}, stratification by the
diagram of initial exponents of the ideal $\Fb(X)$, $b \in X$ 
\cite[Theorem 8.1]{BMAnn}, and a 
$\cC^\infty$ composite function property \cite[\S1.5]{BMAnn}. A uniform
bound for the Chevalley function
measures loss of differentiability in a $\cC^r$ version of 
the composite
function theorem. We use the techniques of \cite{BMAnn} to prove 
Theorems \ref{thm:main1} and \ref{thm:main2} here.

Wang \cite[Theorem~1.1]{Wang} used \cite[Theorem~1.2]{Iz1} to prove that the
Chevalley function associated to a regular proper real-analytic mapping
$\vp\colon M\to\R^n$ has a uniform linear bound  
if and only if $X=\varphi(M)$ has a \emph{uniform
linear product estimate}; i.e., for every compact $K \subset X$, there exist
$\al_K, \be_K \in \N$ such that, for all $b \in K$ and $F,G \in \Oh_b$,
$$
\nu_{X_i,b}(F\cdot G)\ \leq\ \al_K (\nu_{X_i,b}(F)+\nu_{X_i,b}(G)) + \be_K\ ,
$$
where $X_b=\bigcup_iX_i$ is a decomposition of the germ $X_b$ into finitely
many irreducible subanalytic components. We therefore obtain the following
from Theorem \ref{thm:main1}:

\begin{theorem}
\label{thm:prod-estimate}
A closed Nash subanalytic subset of $\R^n$ admits a uniform linear product
estimate.
\end{theorem}

\begin{remark}
\label{artinrees}
The Artin-Rees lemma can be viewed as a version of Chevalley's lemma
for linear relations over a Noetherian ring $R$: Suppose that
$\Psi:\ E \to G$
is a homomorphism of finitely-generated modules over 
$R$, and let $F \subset G$ denote the 
image of $\Psi$. Let $\mi$ be a maximal ideal of $R$.
Then $F \cap \mi^l G \subset \mi^k F$ if and only if
$\Psi^{-1}(\mi^l G) \subset \Ker \Psi + \mi^k E$. The Artin-Rees lemma
says that there exists $\be \in \N$ such that
$F \cap \mi^{k+\be} G = \mi^k (F \cap \mi^{\be} G)$, for all $k$.
In particular, there is always a \emph{linear Artin-Rees exponent}
$l(k) = k + \be$. Uniform versions of the Artin-Rees lemma were proved
in \cite[Theorem 7.4]{BM3}, \cite{DO}, \cite{Hu}. 
A uniform Artin-Rees exponent for a homomorphism of $\cO_M$-modules,
where $M$ is a real-analytic manifold, measures loss of differentiability
in Malgrange division, in the same way that a uniform bound for the Chevalley 
function relates to composite differentiable functions. (See \cite{BM3}.)
\end{remark}
\smallskip

\section{Techniques}
\label{sec:techniques}

\subsection{Linear algebra lemma}
Let $R$ denote a commutative ring with identity, and let $E$ and $F$ be 
$R$-modules. If $B\in\Hom_R(E,F)$ and $r\in\N$, $r\geq 1$, we define
\[
\ad^rB\ \in\ \Hom_R\left(F,
\, \Hom_R\left({\bigwedge}^rE,{\bigwedge}^{r+1}F\right)\right)
\]
by the formula
\[
(\ad^rB)(\omega)(\eta_1\wedge\dots\wedge\eta_r)\ =
\ \omega\wedge B\eta_1\wedge\dots\wedge B\eta_r\ ,
\]
where $\omega\in F$ and $\eta_1,\dots,\eta_r\in E$. 
($\ad^0B:=\id_F$, the identity mapping of $F$.)
Clearly, if $r>\rk B$ then $\ad^rB=0$, and if $r=\rk B$ then 
$\ad^r B\cdot B=0$. ($\rk B$ means the smallest $r$ such that 
$\bigwedge^sB=0$ for all $s>r$.) If $R$ is a field, then 
$\rk B=\dim\Imm B$, so we get:

\begin{lemma}[{\cite[\S6]{BM1}}]
\label{lem:alg-lemma}
Let $E$ and $F$ be finite-dimensional vector spaces over a field $\K$. 
If $B\colon E\to F$ is a linear transformation and $r=\rk B$, then
\[
\Imm B\ =\ \Ker\ad^r B\ .
\]
In particular, if $A$ is another linear transformation with target $F$, 
then $A\xi+B\eta=0$ (for some $\eta$) if and only if 
$\xi\in\Ker\ad^r B\cdot A$.
\end{lemma}
\smallskip

\subsection{The diagram of initial exponents}
\label{sec:diagram}
Let $A$ be a commutative ring with identity. Consider the total ordering
of $\N^n$ given by the lexicographic ordering of $(n+1)$-tuples 
$(|\beta|,\beta_1,\dots,\beta_n)$, where 
$\beta=(\beta_1,\dots,\beta_n)\in\N^n$ and 
$|\beta|=\beta_1+\dots+\beta_n$. For any formal power series 
$F(Y)=\sum_{\beta\in\N^n}F_{\beta}Y^{\beta}\in 
A\lbr Y\rbr =A\lbr Y_1,\dots,Y_n\rbr $, we define the \emph{support} 
$\supp F:=\{\beta\in\N^n\colon F_{\beta}\neq0\}$ and the 
\emph{initial exponent} $\exp F:=\min\,\supp F$. ($\exp F:=\infty$ 
if $F=0$.) 

Let $I$ be an ideal in $A\lbr Y\rbr $. The \emph{diagram of initial 
exponents} of $I$ is defined as
\[
\Ni(I)\ :=\ \{\exp F\colon F\in I\setminus\{0\}\}\ .
\]
Clearly, $\Ni(I)+\N^n=\Ni(I)$.

Suppose that $A$ is a field $\K$. Then, by the formal division theorem 
of Hironaka \cite{Hi1} (see \cite[Theorem~6.2]{BM3}),
\begin{equation}
\label{eqn:div}
\K\lbr Y\rbr \ =\ I\oplus\K\lbr Y\rbr^{\Ni(I)},
\end{equation}
where $\K\lbr Y\rbr^{\Ni}$ is defined as $\{F\in\K\lbr Y \rbr \colon
\, \supp F\subset\N^n\setminus\Ni\}$, for any $\Ni\in\N^n$ such that
$\Ni+\N^n=\Ni$.
\smallskip

\subsection{Fibred product}
\label{sec:fibre}
Let $M$ denote an analytic manifold over $\K$, and let $s\in\N$, $s\geq1$.
Let $\vp\colon\ M \to N$ be an analytic mapping. 
We denote by $\Ms$ the $s$-fold \emph{fibred
product} of $M$ with itself \emph{over} $N$; i.e.,
\[
\Ms\ :=\ \{\ua=(a^1,\dots,a^s)\in M^s\colon\ \vp(a^1)=\dots=\vp(a^s)\}\ ;
\]
$\Ms$ is a closed analytic subset of $M^s$. There is a natural mapping 
$\uvp=\uvp^s\colon$ $\Ms\to N$ given by $\uvp(\ua)=\vp(a^1)$; 
i.e., for each $i=1,\ldots,s$, $\uvp=\vp\circ\rho^i$, where $\rho^i\colon 
\ \Ms \ni (x^1,\dots,x^s) \mapsto x^i\in M$.

Suppose that $\K = \R$.
Let $E$ be a closed subanalytic subset of $M$, and let 
$\varphi\colon\ E\to\R^n$ be a continuous subanalytic mapping. 
Then the fibred product $\Es$ 
is a closed subanalytic subset of $M^s$, and the canonical mapping 
$\uvp=\uvp^s\colon\ \Es\to\R^n$ is subanalytic.

Let $\Esr$ denote the subset of $\Es$ consisting of points 
$\ux=(x^1,\dots,x^s)\in \Es$ such that each $x^i$ lies in a distinct 
connected component of the fibre $\vp^{-1}(\uvp(\ux))$.
If $\vp$ is proper, then $\Esr$ is a subanalytic subset of $M^s$
\cite[\S7]{BMAnn}.
\smallskip

\subsection{Jets}
\label{sec:jets}
Let $N$ denote an analytic manifold (over $\K = \R$ or $\C$), and let $b \in N$.
Let $l\in\N$ and let $J^l(b)$ denote $\Oh_b/\mih_b^{l+1}$. If 
$F\in\Oh_b$, then $J^lF(b)$ denotes the image of $F$ in $J^l(b)$. 
Let $M$ be an analytic manifold, and let $\varphi
\colon M\to N$ be an analytic mapping. 
If $a\in\varphi^{-1}(b)$, then the homomorphism $\vph^*_a\colon
\Oh_b\to\Oh_a$ induces a linear transformation 
$J^l\varphi(a)\colon J^l(b)\to J^l(a)$.

Suppose that $N=\K^n$. Let $y=(y_1,\dots,y_n)$ denote the affine 
coordinates of $\K^n$. Taylor series expansion induces an identification 
of $\Oh_b$ with the ring of formal power series $\K\lbr y-b\rbr =
\K\lbr y_1-b_1,\dots,y_n-b_n\rbr$ (we write $F(y) = 
\sum_{\be \in \N^n} F_{\be}(y-b)^\be$), and hence an identification of 
$J^l(b)$ with $\K^q$, $q=\binom{n+l}{l}$, with respect to which 
$J^lF(b)=(D^{\beta}F(b))_{|\beta|\leq l}$, where $D^{\beta}$ denotes 
$1/\be!$ times the formal derivative of order $\beta \in \N$. 

Using a system of coordinates $x=(x_1,\dots,x_m)$ for $M$ in a 
neighbourhood of $a$, we can identify $J^l(a)$ with $\K^p$, 
$p=\binom{m+l}{l}$. Then
\begin{equation*}
J^l\varphi(a)\colon\ (F_{\beta})_{|\beta|\leq l}\ \mapsto
\ \left((\vph^*_a(F))_{\alpha}\right)_{|\alpha|\leq l}\ =
\ \left(\sum_{|\beta|\leq l}F_{\beta}
L^{\beta}_{\alpha}(a)\right)_{|\alpha|\leq l}\ ,
\end{equation*}
where $L^{\beta}_{\alpha}(a)=
(\partial^{|\alpha|}\varphi^{\beta}/\partial x^{\alpha})(a)/\alpha!$ 
and $\varphi^{\beta}=\varphi_1^{\beta_1}\dots\varphi_n^{\beta_n}$
($\vp = (\vp_1, \dots, \vp_n)$).

Set $J^l_b:=J^l(b)\otimes_{\K}\Oh_b=\bigoplus_{|\beta|\leq l}\K\lbr y-b\rbr$. 
We put $J^l_bF(y):=(D^{\beta}F(y))_{|\beta|\leq l}\in J^l_b$. 
(Evaluating at $b$ transforms $J^l_bF$ to $J^lF(b)$.) 
The ring homomorphism $\vph^*_a\colon\Oh_b\to\Oh_a$ induces a homomorphism 
of $\K\lbr x-a\rbr$-modules,
\begin{equation*}
J^l_a\varphi\colon\quad
\begin{aligned}[t]
&J^l(b)\otimes_{\K}\Oh_a\\ &\qquad\,\,\,\|\\ &\bigoplus_{|\beta|\leq l}
\K\lbr x-a\rbr
\end{aligned}
\quad\longrightarrow\quad
\begin{aligned}[t]
&J^l(a)\otimes_{\K}\Oh_a\\ &\qquad\,\,\,\|\\ &\bigoplus_{|\alpha|\leq l}
\K\lbr x-a\rbr
\end{aligned}
\end{equation*}
such that, if $F\in\Oh_b$, then
\[
J^l_a\varphi\left((\vph^*_a(D^{\beta}F))_{|\beta|\leq l}\right)\ =
\ (D^{\alpha}(\vph^*_a(F)))_{|\alpha|\leq l}.
\]
By evaluation at $a$, $J^l_a\varphi$ induces 
$J^l\varphi(a)\colon J^l(b)\to J^l(a)$. $J^l_a\varphi$ identifies with 
the matrix (with rows indexed by $\alpha\in\N^m$, $|\alpha|\leq l$, 
and columns indexed by $\beta\in\N^n$, $|\beta|\leq l$) whose entries 
are the Taylor expansions at $a$ of the 
$D^{\alpha}\varphi^{\beta}=
(\partial^{|\alpha|}\varphi^{\beta}/\partial x^{\alpha})/\alpha!$, 
$|\alpha|\leq l$, $|\beta|\leq l$.

Let $\una=(a^1,\dots,a^s)\in \Ms$ and let $b=\uvp(\una)$. 
For each $i=1,\dots,s$, the homomorphism $J^l_b = 
J^l(b)\otimes_{\K}\Oh_b\to J^l(a^i)\otimes_{\K}\Oh_{a^i} = J^l_{a^i}$ 
over $\vph^*_{a^i}$, as defined above (using a coordinate system
$x^i = (x^i_1,\dots, x^i_m)$ for $M$ in a neighbourhood of $a^i$), 
followed by the canonical homomorphism 
$J^l(a^i)\otimes_{\K}\Oh_{a^i}\to J^l(a^i)\otimes_{\K}\Oh_{\Ms,\una}$ 
over $(\hat{\rho^i})^*_{\una}\colon\Oh_{a^i}\to\Oh_{\Ms,\una}$, 
induces an $\Oh_{\Ms,\una}$-homomorphism 
$J^l(b)\otimes_{\K}\Oh_{\Ms,\una}\to J^l(a^i)\otimes_{\K}\Oh_{\Ms,\una}$. 
We thus obtain an $\Oh_{\Ms,\una}$-homomorphism
\begin{equation*}
J^l_{\una}\varphi\colon\quad
\begin{aligned}[t]
&J^l(b)\otimes_{\K}\Oh_{\Ms,\una}\\[2ex] &\qquad\,\,\,\|\\[1.5ex] 
&\bigoplus_{|\beta|\leq l}\Oh_{\Ms,\una}
\end{aligned}
\quad\longrightarrow\quad
\begin{aligned}[t]
&\bigoplus_{i=1}^sJ^l(a^i)\otimes_{\K}\Oh_{\Ms,\una}\\ &\qquad\,\,\,\|\\ 
&\bigoplus_{i=1}^s\bigoplus_{|\alpha|\leq l}\Oh_{\Ms,\una} \quad.
\end{aligned}
\end{equation*}

For any (germ at $\una$ of an) analytic subspace $L$ of $\Ms$, 
we also write
\begin{equation}
\label{eqn:J-l-phi-L}
J^l_{\una}\varphi\colon\ J^l(b)\otimes_{\K}\Oh_{L,\una}
\ \to\ \bigoplus_{i=1}^s J^l(a^i)\otimes_{\K}\Oh_{L,\una}
\end{equation}
for the induced $\Oh_{L,\una}$-homomorphism. Evaluation at 
$\una$ transforms $J^l_{\una}\varphi$ to
\begin{equation}
\label{eqn:J-l-phi-una}
J^l\varphi(\una)=(J^l\varphi(a^1),\dots,J^l\varphi(a^s))
\colon\ J^l(b)\to\bigoplus_{i=1}^sJ^l(a^i).
\end{equation}

\section{Ideals of relations and Chevalley functions}
\label{sec:relations}

Let $M$ denote an analytic manifold (over $\K=\R$ or $\C$), and let 
$\varphi=(\varphi_1,\dots,\varphi_n)\colon M \to \K^n$ be an analytic mapping.
If $a\in M$, let $\mathcal{R}_a$ denote the ideal of formal relations
$\Ker\vph^*_a$. 

\begin{remark}
\label{rem:Ra-const-on-comps}
$\mathcal{R}_a$ is constant on connected components of the fibres of $\varphi$
\cite[Lemma~5.1]{BMAnn}.
\end{remark}

Let $s$ be a positive integer, and let $\una=(a^1,\dots,a^s)\in \Ms$. Put
\begin{equation}
\label{eqn:R-una}
\Runa\ :=\ \bigcap_{i=1}^s{\mathcal R}_{a^i}\ =
\ \bigcap_{i=1}^s\Ker\vph^*_{a^i}\ \subset\ \Oh_{\uvp(\una)}\ .
\end{equation}
If $k\in\N$, we also write
\[
\mathcal{R}^k(\una)\ :=
\ \frac{\Runa+\mih^{k+1}_{\uvp(\una)}}{\mih^{k+1}_{\uvp(\una)}}\ 
\subset\ J^k(\uvp(\una))\ .
\]
If $b\in\K^n$, let $\pi^k(b)\colon\Oh_b\to J^k(b)$ denote the canonical 
projection. For $l\geq k$, let $\pi^{lk}(b)\colon J^l(b)\to J^k(b)$ be 
the projection. Set
\[
E^l(\una)\ :=\ \Ker J^l\varphi(\una),\ \ \mathrm{and}
\ \ E^{lk}(\una)\ :=\ \pi^{lk}(\uvp(\una)).E^l(\una)\ .
\]
\smallskip

\subsection{Chevalley's lemma}

\begin{lemma}[{\cite[Lemma~8.2.2]{BM3}}; cf. {\cite[\S~II, Lemma~7]{Chev}}]
\label{lem:Chevalley}
Let $\una\in \Ms$, $\una=(a^1,\dots,a^s)$. For all $k\in\N$, there exists 
$l\in\N$ such that $\mathcal{R}^k(\una)=E^{lk}(\una)$; i.e., such that if 
$F\in\Oh_{\uvp(\una)}$ and $\vph^*_{a^i}(F)\in\mih^{l+1}_{a^i}$, 
$i=1,\dots,s$, then $F\in\Runa+\mih^{k+1}_{\uvp(\una)}$.
\end{lemma}

We write $l(\una,k)=l_{\varphi^*}(\una,k)$ for the least $l$ satisfying 
the conclusion of the lemma.

\begin{proof}[Proof of Lemma \ref{lem:Chevalley}]
If $k\leq l_1\leq l_2$, then
\[
\mathcal{R}^k(\una)\ \subset\ E^{l_2,k}(\una)\ \subset\ E^{l_1,k}(\una)\ ,
\]
and the projection $\pi^{l_2,l_1}(\uvp(\una))$ maps 
$\bigcap_{l\geq l_2}E^{ll_2}(\una)$ onto $\bigcap_{l\geq l_1}E^{ll_1}(\una)$. 
It follows that $\mathcal{R}^k(\una)=\bigcap_{l\geq k}E^{lk}(\una)$. 
Since $\dim J^k(\uvp(\una))<\infty$, there exists $l\in\N$ such 
that $\mathcal{R}^k(\una)=E^{lk}(\una)$.
\end{proof}
\smallskip

\subsection{Generic Chevalley function}
Let $\ua \in \Ms$ and $k \in \N$.
Set
\[
H_{\ua}(k)\ :=\ \dim_{\K}\frac{J^k(\uvp(\una))}{\mathcal{R}^k(\ua)}\ ,
\qquad
d^{lk}(\una)\ :=\ \dim_{\K}\frac{J^k(\uvp(\una))}{E^{lk}(\una)}\ ,
\ \mathrm{if} 
\ l \geq k\ 
\]
($H_{\ua}$ is the \emph{Hilbert-Samuel function} of $\Oh_{\uvp(\ua)}/\Runa$).

\begin{remark}
\label{rem:d-less-than-H}
$d^{lk}(\una)\leq H_{\una}(k)$ since $\mathcal{R}^k(\una)\subset E^{lk}(\una)$. $\mathcal{R}^k(\una)=E^{lk}(\una)$ (and $d^{lk}(\una)= H_{\una}(k)$) if and only if $l\geq l(\una,k)$.
\end{remark}

\begin{lemma}[{\cite[Lemma~8.3.3]{BM3}}]
\label{lem:8.3.3}
Let $L$ be a \emph{subanalytic leaf} in $\Ms$ (i.e., a connected subanalytic
subset of $\Ms$ which is an analytic submanifold of $M^s$; see Remark 
\ref{rem:real}). Then there is a
residual subset $D$ of $L$ such that, if
$\ua,\ua'\in D$, then $H_{\ua}(k)=H_{\ua'}(k)$ and $l(\ua,k)=l(\ua',k)$,
for all $k\in\N$.
\end{lemma}

\begin{definition}
\label{def:gen-estimate}
We define the \emph{generic Chevalley function} of $L$ as 
$l(L,k):=l(\ua,k)$ ($k\in\N$), where $\ua\in D$.
\end{definition}

\begin{proof}[Proof of Lemma \ref{lem:8.3.3}]
For $\ua \in \Ms$ and $l \geq k$,
write $J^l\varphi(\una)$ (\ref{eqn:J-l-phi-una}) (using local coordinates
for $M^s$ as in \S2.4, in a neighbourhood of a point of $\overline{L}$)
as a block matrix
\[
\begin{aligned}
J^l\varphi(\una) \ \ &= \ \ (S^{lk}(\una),T^{lk}(\una))\\
 &= \ \ \left( \begin{array}{cc} J^k\varphi(\una) 
& 0\\ * & * \end{array} \right)
\end{aligned}
\]
corresponding to the decomposition of vectors 
$\xi=(\xi_{\beta})_{\beta\in\N^n, |\beta|\leq l}$ in the source as 
$\xi=(\xi^k,\zeta^{lk})$, where $\xi^k=(\xi_{\beta})_{|\beta|\leq k}$ 
and $\zeta^{lk}=(\xi_{\beta})_{k<|\beta|\leq l}$. Then
\[
E^{lk}(\ua)\ = \ \{\eta=(\eta_{\beta})_{|\beta|\leq k}\colon
\ S^{lk}(\ua)\cdot\eta\in\Imm T^{lk}(\ua)\}\ .
\]
Thus, by Lemma~\ref{lem:alg-lemma}
\[
E^{lk}(\ua)\ = \ \Ker\Theta^{lk}(\una),\ \ \mathrm{and}
\ \ d^{lk}(\ua)\ = \ \rk \Theta^{lk}(\ua)\ ,
\]
where
\[
\Theta^{lk}(\una)\ := \ \ad^{r^{lk}(\ua)}T^{lk}(\ua)\cdot S^{lk}(\ua)\ ,
\quad r^{lk}(\una)\ := \ \rk T^{lk}(\ua)\ .
\]

Set
\[
r^{lk}(L)\ := \ \max_{\una\in L}r^{lk}(\una),\ \ \mathrm{and}
\ \ d^{lk}_L(\una)\ := \ \mathrm{rk}\,\Theta^{lk}_L(\una),\ \ \ua \in L\ ,
\]
where
\[
\Theta^{lk}_L(\una)\ :=\ \mathrm{ad}^{r^{lk}(L)}T^{lk}(\una)\cdot S^{lk}(\una)
\]
(so that $\Theta^{lk}_L(\una)=0$ if $r^{lk}(\una)<r^{lk}(L)$). Let $Y^{lk} := 
\{\una\in L\colon\ r^{lk}(\una)<r^{lk}(L)\}$. Set
\[
d^{lk}(L)\ :=\ \max_{\ua\in L}d^{lk}_L(\ua)\ .
\]
Clearly, $d^{lk}_L(\ua)=0$ if $\ua\in Y^{lk}$, and 
$d^{lk}_L(\ua)=d^{lk}(\ua)$ if $\ua\in L\setminus Y^{lk}$. Also set
\[
Z^{lk}\ :=\ 
Y^{lk}\cup\left\{\una\in L\colon d^{lk}_L(\una)<d^{lk}(L)\right\}\ .
\]
Then $Y^{lk}$ and $Z^{lk}$ are proper closed analytic subsets of $L$. 
For all $\una\in L\setminus Z^{lk}$, $r^{lk}(\una)=r^{lk}(L)$ and 
$d^{lk}(\una)=d^{lk}_L(\una)=d^{lk}(L)$. Put
\begin{equation}
\label{eqn:D}
D^k\ :=\ L\setminus\bigcup_{l>k}Z^{lk}\ , \quad
D\ :=\ \bigcap_{k\geq1}D^k\ .
\end{equation}
By the Baire Category Theorem, the $D^k$ (and hence also $D$) are 
residual subsets of $L$.

Fix $k\in\N$. If $\ua\in D^k$, then $d^{lk}(\ua)=d^{lk}(L)$, for all $l>k$. 
If, in addition, $l\geq l(\ua,k)$, then $H_{\ua}(k)=d^{lk}(L)$, by 
Remark~\ref{rem:d-less-than-H}. If $\ua,\ua'\in D^k$, then, choosing 
$l\geq l(\ua,k)$ and $\geq l(\ua',k)$, we get $H_{\ua}(k)=H_{\ua'}(k)$. 
For the second assertion of the lemma, suppose that $l\geq l(\ua,k)$. Then 
$H_{\ua'}(k)=H_{\ua}(k)=d^{lk}(\ua)=d^{lk}(L)=d^{lk}(\ua')$, so that 
$l\geq l(\ua',k)$, by Remark~\ref{rem:d-less-than-H}. In the same way, 
$l\geq l(\ua',k)$ implies that $l\geq l(\ua,k)$.
\end{proof}
\smallskip

\subsection{Chevalley function of a subanalytic set}
Let $N$ denote a real-analytic manifold, and let $X$ be a closed subanalytic 
subset of $N$. If $b\in X$, then $\Fb(X)$ or $\mathcal{R}_b\ \subset\ \Oh_b$
denotes the \emph{formal local ideal} of $X$ at $b$, in the sense of the 
following simple lemma:

\begin{lemma}
\label{lem:formal-closure}
Let $b\in X$. The following three definitions of $\Fb(X)$ are equivalent:
\begin{itemize}
\item[(1)] Let $M$ be a real-analytic manifold and let $\vp\colon M\to N$ 
be a proper real-analytic mapping such that $X=\vp(M)$. Then 
$\Fb(X)=\bigcap_{a\in\varphi^{-1}(b)}\ker\vph^*_a$.
\item[(2)] $\Fb(X)=\{F\in\Oh_b\colon\ (F\circ\gamma)(t)\equiv0$ for every 
real-analytic arc $\gamma(t)$ in $X$ such that $\gamma(0)=b\}$.
\item[(3)] $\Fb(X)=\{F\in\Oh_b\colon\ T^k_bF(y)=o(|y-b|^k)$, where $y\in X$, 
for all $k\in\N\}$. Here $T^k_bF(y)$ denotes the Taylor polynomial of 
order $k$ of $F$ at $b$, in any local coordinate system.
\end{itemize}
\end{lemma}

Assume that $N = \R^n$, with coordinates $y=(y_1,\dots,y_n)$.
Let $b\in X$. Recall (\ref{eqn:ord}).

\begin{remark}
\label{rem:nu-less-than-mu}
$\nu_{X,b}(F)\leq\mu_{X,b}(F)$: Suppose that $F\in\mih^l_b+\Fb(X)$; 
say $F=G+H$, where $G\in\mih^l_b$ and $H\in\Fb(X)$. Then 
$|T^l_bG(y)|\leq c|y-b|^l$ and $T^l_bH(y)=o(|y-b|^l)$, $y\in X$, 
by Lemma~\ref{lem:formal-closure}. Hence 
$|T^l_bF(y)|\leq\,\mathrm{const}|y-b|^l$ on $X$.
\end{remark}

\begin{definition}[\emph{Chevalley functions}]
\label{def:Chev-estimates}
Let $b\in X$ and let $k\in\N$. Set
\[
l_X(b,k)\ :=\ \min\{l\in\N\colon\mathrm{\ If\ }F\in\Oh_b\mathrm{\ and\ }
\mu_{X,b}(F) > l,\mathrm{\ then\ }\nu_{X,b}(F) > k\}\ .
\]
Let $\vp\colon M\to N$ be a proper real-analytic mapping such that 
$X=\vp(M)$. Set
\[
\begin{aligned}
l_{\vp^*}(b,k)\ :=\ \min\{l\in\N\colon
&\mathrm{\ If\ }F\in\Oh_b\mathrm{\ and\ }\nu_{M,a}(\vph^*_a(F)) > l\\
&\mathrm{for\ all\ }a\in\vp^{-1}(b),\mathrm{\ then\ }\nu_{X,b}(F) > k\}\ .
\end{aligned}
\]
\end{definition}

\begin{remark}
\label{rem:l-b-less-than-l-a}
Suppose that $b = \uvp(\ua)$, where $\ua=(a^1,\dots,a^s) \in \Ms$, $s \geq 1$.
By Lemma~\ref{lem:Chevalley}, $l_{\vp^*}(\ua,k)<\infty$. If $\ua$ 
includes a point $a^i$ in every connected component of $\vp^{-1}(b)$, 
then $\bigcap_{i=1}^s\Ker\vph^*_{a^i}=\Fb(X)$ (by 
Remark \ref{rem:Ra-const-on-comps} and
Lemma~\ref{lem:formal-closure}), so that 
$l_{\vp^*}(b,k)\leq l_{\vp^*}(\ua,k)$.
\end{remark}

\begin{lemma}[see {\cite[Lemma~6.5]{BMAnn}}]
\label{lem:6.5}
Let $\varphi\colon M\to N$ be a proper real-analytic mapping such that $X=\varphi(M)$. Then
$l_X(b,\cdot)\leq l_{\varphi^*}(b,\cdot)$ for all $b\in X$.
\end{lemma}
\smallskip

\section{Proofs of the main theorems}
\label{sec:uniform-lin-est}

Let $\vp\colon M\to\K^n$ be an analytic mapping from a manifold $M$ (over 
$\K=\R$ or $\C$). Let $s$ be a positive integer. Let $\ua=(a^1,\dots,a^s)
\in \Ms$, and let $b=\uvp(\ua)$. 

\begin{remark}
\label{rem:estimates-supported-outside}
By (\ref{eqn:div}), the Chevalley functions $l_{\varphi^*}(\una,k)$ and 
$l_{\varphi^*}(b,k)$ (Definitions~\ref{def:Chev-estimates}) can be defined 
using power series that are supported outside the diagram of 
initial exponents: Set $\Ni_{\ua}:=\Ni(\Runa)$ and $\Ni_b:=\Ni(\mathcal{R}_b)$
(cf.~\ref{eqn:R-una} and Lemma~\ref{lem:formal-closure}). Then
\[
\begin{aligned}
l_{\vp^*}(\ua,k)\ &=\ \begin{aligned}[t]\min\{l\in\N\colon\ &\mathrm{If\ }
F\in\Oh^{\Ni_{\ua}}_b\mathrm{\ and\ }\vph^*_{a^i}(F)\in\mih^{l+1}_{a^i},
\ i=1,\dots,s,\\ &\mathrm{then\ }F\in\Runa+\mih^{k+1}_b\}\ ,\end{aligned}\\
l_{\vp^*}(b,k)\ &=\ \begin{aligned}[t]\min\{l\in\N\colon\ &\mathrm{If\ }
F\in\Oh^{\Ni_b}_b\mathrm{\ and\ }\vph^*_a(F)\in\mih^{l+1}_a,\mathrm{\ for\ all\ }
a\in\vp^{-1}(b),\\ &\mathrm{then\ }F\in\mathcal{R}_b+\mih^{k+1}_b\}\ .
\end{aligned}
\end{aligned}
\]
(In the latter, we assume that $\vp$ is a proper real-analytic mapping.)
\end{remark}

If $l\in\N$, set $J^l(b)^{\Ni_{\ua}}:=
\{\xi=(\xi_{\beta})_{|\beta|\leq l}\in J^l(b)\colon\xi_{\beta}=0$ 
if $\beta\in\Ni_{\ua}\}$. Consider the linear mapping
\[
\Phi^l(\ua)\colon\ J^l(b)^{\Ni_{\ua}} \to \bigoplus_{i=1}^sJ^l(a^i)
\]
obtained by restriction of $J^l\vp(\ua)\colon J^l(b)\to\bigoplus J^l(a^i)$ 
(\ref{eqn:J-l-phi-una}). Given $k\leq l$, write $\Phi^l(\ua)$ as a 
block matrix
\[
\Phi^l(\ua)\ =\ (A^{lk}(\ua),B^{lk}(\ua))\ ,
\]
where $A^{lk}(\ua)$ is given by the restriction of $\Phi^l(\ua)$ to 
$J^k(b)^{\Ni_{\ua}}$.

\begin{remark}
\label{rem:8.13}
If $\xi\in J^l(b)^{\Ni_{\ua}}$, write $\xi=(\eta,\zeta)$ corresponding to 
this block decomposition. Then 
$l\geq l_{\vp^*}(\ua,k)$ if and only if $A^{lk}(\ua)\eta+B^{lk}(\ua)\zeta=0$ 
implies $\eta=0$ \cite[Lemma~8.13]{BMAnn}.
\end{remark}

\begin{lemma}[(cf. {\cite[Prop.~8.15]{BMAnn}}]
\label{lem1}
Let $s \geq 1$ and consider $\uvp = \uvp^s\colon \Ms \to \R^n$. Let $L$ be
a relatively compact subanalytic leaf in $\Ms$ (cf. Lemma \ref{lem:8.3.3})
such that $\Ni_{\ua}=\Ni(\Runa)$ is constant on $L$. Let $l(k) = l(L,k)$
denote the generic Chevalley function of $L$. Then there exists $p
\in \N$ such that $l_{\vp^*}(\ua,k) \leq l(k) + p$, for all $\ua \in L$
and $k \in \N$.
\end{lemma}

\begin{proof}
Set $\Ni = \Ni_{\ua}$, $\ua \in L$.
We can assume that $\overline{L}$ lies in a coordinate chart for $M^s$
as in \S2.4.
Let $k\in\N$ and let $l=l(k)$. Let $\ua=(a^1,\dots,a^s)\in L$, and set 
$b=\uvp(\ua)$. Consider the linear mapping $\Phi^l(\ua)=
(A^{lk}(\ua),B^{lk}(\ua))\colon J^l(b)^{\Ni}\to\bigoplus_{i=1}^sJ^l(a^i)$ 
as above. The $\Oh_{L,\ua}$-homomorphism $J^l_{\ua}\vp\colon 
J^l(b)\otimes_{\K}\Oh_{L,\ua}\to
\bigoplus_{i=1}^s J^l(a^i)\otimes_{\K}\Oh_{L,\ua}$ (\ref{eqn:J-l-phi-L})
induces an $\Oh_{L,\ua}$-homomorphism
\[
\Phi^l_{\ua}\ =\ (A^{lk}_{\ua},B^{lk}_{\ua})\colon
\ J^l(b)^{\Ni}\otimes_{\K}\Oh_{L,\ua}\ \to
\ \bigoplus_{i=1}^s J^l(a^i)\otimes_{\K}\Oh_{L,\ua}\ ;
\]
evaluating at $\ua$ transforms $\Phi^l_{\ua}$ to $\Phi^l(\ua)\ =
\ (A^{lk}(\ua),B^{lk}(\ua))$.

Let $r=\rk B^{lk}_{\ua}=\mathrm{\ generic\ rank\ of\ }B^{lk}(\ux)$, 
$\ux \in L$. Let $\Theta_{\ua}=\ad^r B^{lk}_{\ua}\cdot A^{lk}_{\ua}$. 
Then $\Ker\Theta_{\ua}=0$ (i.e., $\Ker\Theta(\ux)=0$ generically on $L$, 
where $\Theta(\ux)=\ad^r B^{lk}(\ux)\cdot A^{lk}(\ux)$, by 
Remark~\ref{rem:8.13}). Let $d=\rk \Theta_{\ua}$. Then there is a 
nonzero minor $\delta_{\ua}\in\cO_{L,\ua}$ of $\Theta_{\ua}$ of order $d$; 
$\delta_{\ua}$ is induced by a minor $\delta(\ux)$ of order $d$ of 
$\Theta(\ux)$, $\ux\in L$, such that $\delta(\ux)\neq0$ on a residual 
subset of $L$. Since $\delta$ is a restriction to $L$ of an analytic 
function defined in a neighbourhood of $\overline{L}$, the order of 
$\delta_{\ux}$, $\ux\in L$, is bounded on $L$; say, $\delta_{\ux}\leq p$.

We claim that $l_{\vp^*}(\ua,k)\leq l(k)+p$ for all $\ua\in L$: 
Let $\ua=(a^1,\dots,a^s)\in L$, and let $b=\uvp(\ua)$. Let $l=l(k)$ and 
$l'=l+p$. Suppose that $F\in\Oh^{\Ni}_b$ and 
$\vph^*_{a^i}(F)\in\mih^{l'+1}_{a^i}$, $i=1,\dots,s$. Let 
$\hat\xi_{\ua}=(\hat\eta_{\ua},\hat\zeta_{\ua})$ denote the element of 
$J^l(b)^{\Ni}\otimes_{\K}\Oh_{L,\ua}$ induced by 
$J^l_bF\in J^l(b)\otimes_{\K}\Oh_b$ via the pull-back. Then each component 
of $A^{lk}_{\ua}\hat\eta_{\ua}+B^{lk}_{\ua}\hat\zeta_{\ua}$ belongs to 
$\mih^{l'+1-l}_{L,\ua}$ (as we see by taking formal derivatives of 
order $\leq l$ of the $\vph^*_{a^i}(F)$). It follows that each component 
of $\Theta_{\ua}\hat\eta_{\ua}$ and therefore (by Cramer's rule) each 
component of $\delta_{\ua}\cdot\hat\eta_{\ua}$ belongs to 
$\mih^{l'+1-l}_{L,\ua}$. Thus, each component of $\hat\eta_{\ua}$ lies 
in $\mih^{l'+1-l-p}_{L,\ua}=\mih_{L,\ua}$; i.e., 
$\hat\eta_{\ua}(\ua)=0$, so that $F$ vanishes to order $k$ at $b=\uvp(\ua)$.
\end{proof}

\begin{proof}[{\bf Proof of Theorem {\ref{thm:main1}}}]
By \cite[Theorems~A,C]{BM3}, there is a locally finite partition of $M$
into relatively compact subanalytic leaves $L$ such that the diagram
of initial exponents $\Ni_{a}=\Ni(\Ra)$ is constant on each $L$.
Given $L$, let $l(L,k)$ denote the generic Chevalley function. (In
particular, $l(L,k) = l_{\vp^*}(a,k)$, for all $a$ in a residual
subset of $L$.) Since $\vp$ is regular, there exist $\al_L, \gamma_L$
such that $l(L,k) \leq \al_L k + \gamma_L$, for all $k \in \N$ (by
\cite{Iz2}). By Lemma \ref{lem1} (in the case $s=1$), there exists
$p_L \in \N$ such that $l_{\vp^*}(a,k) \leq \al_L k + \gamma_L + p_L$,
for all $a \in L$ and all $k$. The result follows.
\end{proof}

\begin{remark}
\label{rem:real}
In the case $\K=\C$, we define ``subanalytic leaf'' using the underlying
real structure. If $\vp$ is regular, then the diagram $\Ni_{a}$ is, in fact, 
an upper-semicontinuous function of $a$, with respect to the $\K$-analytic
Zariski topology of $M$ (and a natural total ordering of $\{\Ni \in \N^n\colon
\Ni + \N^n = \Ni\}$) \cite[Theorem C]{BM3}, but we do not     
need the more precise result here.
\end{remark}

\begin{lemma}
\label{lem2}
Let $s \geq 1$ and let $\ua = (a^1,\dots, a^n) \in \Ms$. Suppose that 
$\vp$ is regular at $a^1,\dots, a^n$. Then there exist $\al, \be \in \R$
such that $l_{\vp^*}(\ua, k) \leq \al k + \be$, for all $k \in \N$.
\end{lemma}

\begin{proof}
Let $b = \uvp(\ua)$. For each $i = 1,\dots, s$, since $\vp$ is regular
at $a^i$, there exist $\al^i, \be^i$ such that
\begin{equation}
\label{eq1}
l_{\vp^*}(a^i,k)\ \leq\ \al^i k + \be^i, \quad \mathrm{for\ all\ } k\ .
\end{equation}
Of course, $\bigcap_{i=1}^s\Ker\vph^*_{a^i}$ is the kernel of the
homomorphism $\Oh_b \to \bigoplus_{i=1}^s\Oh_b/\ker\vph^*_{a^i}$. By the
Artin-Rees lemma (cf. Remark \ref{artinrees}), there exists $\lambda \in \N$
such that, if $F \in \mih^{k+\lambda}_b + \ker\vph^*_{a^i}$, $i=1,\dots,s$,
then 
\begin{equation}
\label{eq2}
F\ \in\ \mih^k_b + \bigcap_{i=1}^s\Ker\vph^*_{a^i}\ .
\end{equation}

Now let $F \in \Oh_b$ and suppose that $\vph^*_{a^i}(F) \in \mih^{\al^i(\lambda
+k)+\be^i+1}_{a^i}$, $i=1,\dots,s$. Then $F \in \mih^{\lambda +k+1}_b + 
\Ker \vph^*_{a^i}$, $i=1,\dots,s$, by (\ref{eq1}), so that $F \in
\mih^{k+1}_b + \bigcap_{i=1}^s\Ker\vph^*_{a^i}$, by (\ref{eq2}). In other
words, $l_{\vp^*}(\ua, k) \leq \al k + \be$, where $\al = \max \al^i$
and $\be = \lambda \max \al^i + \max \be^i$.
\end{proof}

\begin{proof}[{\bf Proof of Theorem {\ref{thm:main2}}}]
Suppose that $\vp\colon M \to \R^n$ is a real-analytic mapping, where $M$
is compact. Let $X = \vp(M)$. Let $s \geq 1$, $\ua \in \Ms$, $b = \uvp(\ua)$.
If $\ua = (a^1,\dots, a^s)$ includes a point $a^i$ in every connected component
of $\vp^{-1}(b)$, then
\begin{equation}
\label{eq3}
l_X(b, k)\ \leq\ l_{\vp^*}(\ua, k)\ ,
\end{equation}
by Remark \ref{rem:l-b-less-than-l-a} and Lemma \ref{lem:6.5}.

Let $L$ be a relatively compact subanalytic leaf in $\Ms$, such that
$\Ni_{\ua}=\Ni(\Runa)$ is constant on $L$. Suppose that $\vp$ is regular
at $a^i$, for all $\ua = (a^1,\dots, a^s) \in L$ and $i=1,\dots,s$. Let
$l(L,k)$ denote the generic Chevalley function of $L$. By Lemma \ref{lem2},
there exist $\al, \be$ such that $l(L,k) \leq \al k + \be$. Therefore,
by Lemma \ref{lem1}, there exist $\al_L, \be_L$ such that
\begin{equation}
\label{eq4}
l_{\vp^*}(\ua, k)\ \leq\ \al_L k + \be_L, \quad \mathrm{for\ all\ \ } \ua \in L\ .
\end{equation}

To prove the theorem, we can assume that $X$ is compact. Let $\vp$ be a
mapping as above, such that $X = \vp(M)$. 
We consider first the case that $X$ is Nash. Then we can assume that
$\vp$ is regular. Let $s$ denote a bound on the 
number of connected components of a fibre $\vp^{-1}(b)$, for all $b \in X$.
Then there is a finite partition of $\Ms$ into relatively compact subanalytic
leaves $L$, such that $\Ni_{\ua}=\Ni(\Runa)$ is constant on every $L$.
By (\ref{eq3}) and (\ref{eq4}), for each $L$, there exist $\al_L, \be_L$
such that $l_X(b,k) \leq \al_L k + \be_L$, for all $b \in \uvp(L)$ and 
all $k$. Therefore, $l_X(b,k)$  has a uniform linear bound.

Finally, we consider $X$ formally Nash. 
Let $\NR(\vp)\subset M$ denote the set of points at which $\vp$ is not
regular. Then $\NR(\vp)$ is a nowhere-dense closed analytic subset of $M$ 
(\cite[Theorem~1]{P}). For each positive integer $s$, set
\[
\NR(\uvp^s)\ :=
\ \Ms \cap \bigcup_{i=1}^s\{\ua = (a^i,\dots,a^s) \in M^s\colon 
a^i\in \NR(\vp)\}\ ;
\]
then $\NR(\uvp^s)$ is a closed analytic subset of $\Ms$.

If $b\in X$ and $a, a'$ belong to the same connected component of 
$\vp^{-1}(b)$, then $\vp$ is regular at $a$ if and only if $\vp$ is 
regular at $a'$ (cf. Remark~\ref{rem:Ra-const-on-comps}). Let $t$ be a
bound on the number of connected components of a fibre $\vp^{-1}(b)$, for
all $b \in X$. For each $s\leq t$, define $X_s:=\{b\in X\colon 
\vp^{-1}(b)$ has precisely $s$ regular components$\}$ and $Y_s:=
\{b\in X\colon \vp^{-1}(b)$ has at least $s$ regular components$\}$. 
Then $X_s=Y_s\setminus Y_{s+1}$, and
\[
Y_s\ =\ \uvp^s(\Msr \setminus \NR(\uvp^s))\ ;
\]
in particular, all the $X_s$ and $Y_s$ are subanalytic (cf. \S3.2).

The hypothesis of the theorem implies:
\begin{enumerate}
\item $X=\bigcup_{s=1}^tX_s$;
\item If $b\in X_s$ and $\ua \in (\uvp^s)^{-1}(b) \bigcap (\Msr \setminus 
\NR(\uvp^s))$, then $\Runa = \mathcal{R}_b$.
\end{enumerate}
((2) follows from the fact that $\Fb(X)=\Fb(Y_b)$, where $Y_b$ is some 
closed Nash subanalytic subset of $X$, and (1) from the fact that the 
latter condition holds for all $b\in X$.) 

By \cite[Theorem~2]{P}, for each $s$, there is a finite stratification 
$\mathcal{L}_s$ of $\Ms$ compatible with $\NR(\uvp^s)$ such that 
$\Ni_{\ua}=\Ni(\Runa)$ is constant on every stratum $L\subset \Ms 
\setminus \NR(\uvp^s)$, $L \in \mathcal{L}_s$. Clearly,
\[
X_s\ = \bigcup_{\substack{L\in\mathcal{L}_s\\L\subset \Ms \setminus 
\NR(\uvp^s)}} \uvp^s \left(L \cap \Msr \right) \cap X_s\ ;
\]
hence
\[
X\ =\ \bigcup_{s=1}^t \bigcup_{\substack{L\in\mathcal{L}_s\\L\subset \Ms
\setminus \NR(\uvp^s)}} \uvp^s \left(L \cap \Msr \right)\ .
\]
Again by (\ref{eq3}) and (\ref{eq4}), for each $L$, there exist $\al_L, 
\be_L$ such that $l_X(b,k) \leq \al_L k + \be_L$, for all $b \in \uvp(L)$
and all $k$. The result follows.
\end{proof}
\smallskip

\bibliographystyle{amsplain}

\end{document}